\documentclass[a4paper, notitlepage, 12pt]{article}

\usepackage[total={160mm,250mm}, top=25mm, left=25mm, includefoot]{geometry}

\usepackage{ucs}              
\usepackage[utf8x]{inputenc} 	

\usepackage[T1]{fontenc}     
\usepackage{lmodern}         

\usepackage[singlespacing]{setspace}

\usepackage{color}
\usepackage{amsmath,amsthm,amssymb,amscd,amsbsy}
\usepackage{upgreek}

\usepackage{tensor}

\usepackage{wasysym}
\usepackage{float}
\usepackage{tabularx,paralist}
\usepackage{graphicx}
\usepackage[breaklinks]{hyperref}
\usepackage{authblk}

\usepackage{multirow} 
\usepackage{array}
\usepackage{caption}

\usepackage{rcs} \RCS $Revision: 0.1 $ \RCS $Date: 2019/11/22 17:59:39 $ \RCS $Author: hgm $ \RCS $RCSfile: 19_Lin-Param-ovr.tex,v $ \RCS $Id: 19_Lin-Param-ovr.tex,v 0.1 2019/11/22 17:59:39 hgm Exp $

\usepackage[british]{babel}

\usepackage{refdef}
\usepackage{ifontdef}

\newcommand{\ignore}[1]{}

\newcommand{\Lp}{\mrm{L}}

\newcommand{\vsigma}{\varsigma}

\newcommand{\vphi}{\varphi}
\newcommand{\vpi}{\varpi}
\newcommand{\vkappa}{\varkappa}

\newcommand{\vek}[1]{\mathchoice{\displaystyle\boldsymbol{#1}}
{\textstyle\boldsymbol{#1}}{\scriptstyle\boldsymbol{#1}}
{\scriptscriptstyle\boldsymbol{#1}}}

\newcommand{\mat}[1]{\mathchoice{\displaystyle\mathbf{#1}}
{\textstyle\mathbf{#1}}{\scriptstyle\mathbf{#1}}
{\scriptscriptstyle\mathbf{#1}}}

\newcommand{\opb}[1]{\vek{{\mathsf{#1}}}}

\newcommand{\ops}[1]{\mathchoice{\displaystyle\mathsf{#1}}
{\textstyle\mathsf{#1}}{\scriptstyle\mathsf{#1}}
{\scriptscriptstyle\mathsf{#1}}}

\newcommand{\tnb}[1]{\mathchoice{\displaystyle\mathboldsans{#1}}
{\textstyle\mathboldsans{#1}}{\scriptstyle\mathboldsans{#1}}
{\scriptscriptstyle\mathboldsans{#1}}}

\newcommand{\tns}[1]{\mathchoice{\displaystyle\mathsans{#1}}
{\textstyle\mathsans{#1}}{\scriptstyle\mathsans{#1}}
{\scriptscriptstyle\mathsans{#1}}}

\newcommand{\EXP}[1]{\mathbb{E}\left(#1\right)}

\newcommand{\diag}{\mathop{\mathrm{diag}}\nolimits}

\newcommand{\im}{\mathop{\mathrm{im}}\nolimits}
\newcommand{\dom}{\mathop{\mathrm{dom}}\nolimits}

\newcommand{\spn}{\mathop{\mathrm{span}}\nolimits}

\newcommand{\cl}{\mathop{\mathrm{cl}}\nolimits}

\newcommand{\di}{\mathrm{d}}

\newcommand{\ip}[2]{\langle #1, #2 \rangle}
\newcommand{\bkt}[2]{\langle #1 \mid #2 \rangle}

\newcommand{\trpos}{{\ops{T}}}

\definecolor{myred}{rgb}{1, 0.2, 0.2}

\newcommand{\autheadcr}{\authorcr}
\newcommand{\citep}[1]{\cite{#1}}

\newcommand{\authorhgm}{Hermann G. Matthies}
\newcommand{\authorro}{Roger Ohayon}

\newcommand{\affilwire}{Institute of Scientific Computing\autheadcr
                        Technische Universit\"at Braunschweig\autheadcr
                        38092 Braunschweig, Germany\autheadcr
                        e-mail: \ttt{wire@tu-bs.de}}
\newcommand{\affilcnam}{Structural Mechanics and Coupled Systems Laboratory\autheadcr
                       Conservatoire National des Arts et M\'etiers (CNAM)\autheadcr
                      75141 Paris Cedex 03, France}

\newcommand{\thetitle}{Parametric Models Analysed with Linear Maps}

\newcommand{\theauthor}{\authorhgm, \authorro}
\newcommand{\thesubject}{35B30, 37M99, 41A05, 41A45, 41A63, 60G20, 60G60, 65J99, 93A30}
\newcommand{\thekeywords}{parametric models, reduced order models
                          reproducing kernel Hilbert space,
                          correlation, factorisation, spectral decomposition,
			                 representation}

\newcommand{\textdate}{\today}

\begin{document}

\title{\thetitle\thanks{Partly supported by the Deutsche
          Forschungsgemeinschaft (DFG) through SPP 1886 and SFB 880.}}


\author[$\dag$]{\authorhgm}
\author[$\ddag$]{\authorro}


\affil[$\dag$]{\affilwire}
\affil[$\ddag$]{\affilcnam}

\date{\textdate}


\maketitle

%

\begin{abstract}
Parametric entities appear in many contexts, be it in optimisation, control,
modelling of random quantities, or uncertainty quantification.  These are all
fields where reduced order models (ROMs) have a place to alleviate the
computational burden.
Assuming that the parametric entity takes values in a linear space, we show
how is is associated to a linear map or operator.  This
provides a general point of view on how to consider and analyse different
representations of such entities.
Analysis of the associated linear map in turn connects such representations
with reproducing kernel Hilbert spaces and affine- / linear-representations in 
terms of tensor products.  A generalised correlation operator is defined
through the associated linear map, and its spectral analysis helps to shed 
light on the approximation properties of ROMs.
This point of view thus unifies many such representations under a functional
analytic roof, leading to a deeper understanding and
making them available for appropriate analysis.

\vspace{5mm}

\vspace{5mm}
{\noindent\textbf{AMS Classification:} \thesubject}

\end{abstract}

%
%
%
%









%

\section{Introduction}  \label{S:intro}
Many mathematical and computational models depend on parameters.  These may be
quantities which have to be optimised during a design, or controlled in a real-time
setting, or these parameters may be uncertain and represent uncertainty
present in the model.  Such parameter dependent models are usually specified in such a way
that an input to the model, e.g.\ a process or a field, depends on these parameters.
In an analogous fashion, the output or the ``state'' of the model will depend on
those parameters.  Any of these entities may be called a parametric model.
To make things a bit more specific, we look at an example:
Consider the parametric entities in the following equation:
\begin{equation}  \label{eq:prim-ex}
    A(\mu; u) = f(\mu).
\end{equation}
Here $A(\zeta(\mu);\cdot):\C{V}\to\C{V}$ is a possibly nonlinear opertor from the
Hilbert space $\C{U}$ into itself, dependent on $\zeta(\mu)\in\C{Z}$ --- a vector
in another Hilbert space $\C{Z}$ used to specify the system --- $u\in\C{V}$ is the
state of the system described by $A$, whereas $f(\mu)\in \C{V}$ is the
excitation resp.\ action on the system.  The \emph{parameters} $\mu\in\C{M}$
are elements of some admissible parameter set $\C{M}$.  
Here $f(\mu)$ and $\zeta(\mu)$ are two examples of such
parametric entities; and as the whole equation depends on $\mu$, we assume that
for each $\mu\in\C{M}$ the system \refeq{eq:prim-ex} will be well-posed and allow
for the state $u(\mu)$ to also be a unique function of the parameters --- another
example of a parametric entity.

When one has to do computations with a system such as \refeq{eq:prim-ex}, one
needs computational \emph{representations} of the parametric entities such as
the ``inputs'' $f(\mu), \zeta(\mu)$, and also the to be determined \emph{state} $u(\mu)$,
the ``output''.  Let us denote any of such generic entities as $r(\mu)$; then one
seeks a computational expression to compute $r(\mu)$ for any given parameter $\mu\in\C{M}$.
The first question which has to be addressed is how to choose ``good co-ordinates''
on the parameter set $\C{M}$.  With this we mean scalar functions $\xi_m:\C{M}\to\C{R}$,
so that the collection and specification of all $\{\xi_m(\mu)\}_{m=1,\dots,M}$
will on one hand specify the particular $\mu\in\C{M}$ as regards the system 
\refeq{eq:prim-ex}, and on the other hand be a computational handle for the
parametric entities $r(\mu)$, which now can be expressed as $r(\xi_1,\dots,\xi_M)$.
Often the parameter set is already given as 
$\C{M}\subseteq \D{R}^d$, so that $\mu = [\eta_1,\dots,\eta_d] \in \D{R}^d$ are
directly given co-ordinates, and the co-ordinate functions $\eta_k$ may directly
serve as co-ordinates.  But often, and not only, but especially, when $d\in\D{N}$
is a large number, it may be advisable to choose other co-ordinates $\xi_m$, 
which should be free of possible constraints and be as ``independent'' as possible.
This is usually part of finding a good computational representation for $r(\mu)$,
and will be addressed as part of our analysis.  One may term this as a 
\emph{re-parametrisation} of the problem.

The second question to be addressed is the actual number of degrees-of-freedom
needed to describe the behaviour of the system \refeq{eq:prim-ex} through
some finite-dimensional approximation or discretisation.  Often the initial
process of discretisation produces a first approximation with a large number
of degrees-of-freedom; this initial computational model is often referred to
as a full-scale or high-fidelity model.
For many computational purposes it is necessary to reduce the number of 
degrees-of-freedom in the computational model in order to be able to
carry out the computations involved in an optimisation
or uncertainty quantification in a acceptable amount of time; such computational
models are then termed reduced order models (ROMs).  If the 
high-fidelity model is a parametric model, the same is required from
the ROM $r_a(\mu)\approx r(\mu)$.

The question of how to produce ROMs for specific kinds of systems like \refeq{eq:prim-ex}
is an important one, and is the subject of extensive current research.  For the general
subject of model order reduction there is an excellent collection of recent work in
\citep{Quarteroni2014} and survey in \citep{chinestaWillcox2017}, as well as an introductory
text in \citep{Quarteroni2015}; see also \citep{FickMadayPatera2017, LamWillcox2018} for
important contributions.  Besides these general considerations, in the present case
parametrised ROMs are of particular interest.  The general survey 
\citep{BennWilcox-paramROM2015} covers the literature up to 2015 very well, 
as well as the later one \citep{ChenSchwab2017}, which is concerned mainly with
uncertainty quantification.  Excellent collections on the topic of parametrised ROMs
are contained in  \citep{MoRePaS2015} and \citep{BennerWillcox2017}.  
A recent systematic monograph is \citep{HesthavenRozza2016}, and important recent
contributions are e.g.\ \citep{BuffaMadayPatera2012, VenturiRozza2019, ZancanaroRozza2019}.
Machine learning and so-called data-driven procedures have also been used in
this context, see the recent contributions in
\citep{HijaziRozza2019, Quarteroni2019, RaissiKarniadakis2019, Schwab2019, SoizeGhanem2019},
but this is at the very beginning.

Here a particular point of view will be taken for the analysis --- not to
be found in the recent literature just surveyed --- namely the
identification of a parametric entity with a linear mapping defined on the
dual space, which is introduced in \refS{parametric}.  
This idea has been around for a long time, and has surfaced
mostly when the ``strong'' notion of a concept has to be replaced by a ``weaker''
one.  In this sense one may see the present point of view as a generalisation
of the view of distributions of generalised functions as linear mappings
\citep{gelfand64-vol1, gelfand64-vol2}.  They were used to define weak notions of
random quantities \citep{gelfand64-vol4}, and some of the present ideas are
also contained in \citep{kreeSoize86}.  In some sense these ideas are already
contained in \citep{Karhunen1947} --- see also the English translation 
\citep{Karhunen1947-e} --- and may most probably be found even earlier.
The reason on why to approach the subject in this way is that for linear operators
there is a host of methods which can be used for their analysis, and it puts
all such parametric entities under one ``roof''.

Here we want to explain the basic abstract framework and how it applies to ROMs.
This present work is a continuation of \citep{boulder:2011} and 
\citep{hgmRO-1-2018, hgmRO-2-2018, hgm-3-2018}.  The general theory was
shown in \citep{hgmRO-1-2018}, and here the purpose is primarily to give
an introduction into this kind of analysis, which draws strongly on the spectral
analysis of self-adjoint operators (e.g.\ \citep{gelfand64-vol3, 
gelfand64-vol4, DautrayLions3}), and an overview on how to use it in the analysis
of ROMs.  This is the topic of \refS{correlat}.
Coupled systems and their ROMs are the focus of \citep{hgmRO-2-2018},
and \citep{hgm-3-2018} is a short note on how this is used for random fields
and processes.  In the \refS{xmpls} some examples of such refinements of
the basic concept are given.

As will be seen, it is very natural to deal with tensor products in
this topic of parametrised entities.  In the form of the 
\emph{proper generalised decomposition} (PGD) this idea has been
explained and used in \citep{chinestaPL2011, AmmarChinestaFalco2010, Falco:2012,
chinestaBook, chinestaWillcox2017}.
The topic of tensor approximations \citep{Hackbusch_tensor} turns
out to be particularly relevant here, and recently new connections between
such approximations and machine learning with deep networks have been
exposed \citep{CohenSha2016, KhrulkovEtal2018}.  In \refS{concl} we
conclude with a recapitulation of the main ideas.


%
%
%
%
%
%



\section{Parametric models and linear maps} \label{S:parametric}
This is a gentle introduction and short recap of the developments in 
\citep{hgmRO-1-2018, hgmRO-2-2018, hgm-3-2018},
where the interested reader may find more detail.
To start, and to take a simple motivating example, one could think of a
scalar function $r(x,\mu)$, defined on some set $\C{X}$, which depends on some
parameters in a set $\C{M}$ --- in other words a parametric function.
In what follows, this function will be viewed as a mapping
\begin{equation}   \label{eq:def-r-map}
   r:\C{M} \to \C{U},
\end{equation}
so that for each value of $\mu\in\C{M}$ the function $w := r(\cdot,\mu)\in\C{U}$
is a scalar function $w(x)$ defined on the set $\C{X}$.

To simplify further and make everything finite-dimensional, assume that we are
only interested in four positions in $\C{X}$, namely $x_1, x_2, x_3, x_4\in\C{X}$,
or, alternatively and even simpler, that
$\C{X}=\{x_1, x_2, x_3, x_4\}$ has only four elements, and finally
for the sake of simplicity, that the parameter
set has only three elements $\C{M} = \{ \mu_1, \mu_2, \mu_3 \}$.
Then one can arrange all the possible values of $r(x,\mu)$ with the abbreviation
$r_{i,j} = r(x_i,\mu_j), (i=1,\dots,3,j=1,\dots,4)$ in the following matrix:
\[ \vek{R}  =  [r_{i,j}]_{i=1,\dots,3,j=1,\dots,4} \in \D{R}^{3 \times 4}.
\]
It is obvious that knowing the function $r(x,\mu)$ is equivalent with knowing
the matrix $\vek{R}$.  As a matrix $\vek{R}$ obviously corresponds to a linear
mapping from $\C{U} =\D{R}^4$ to $\C{R} = \D{R}^3$, and one has for any 
$\vek{u} = [u(x_1), u(x_2), u(x_3), u(x_4)]^\trpos = [u_1, u_2, u_3, u_4]^\trpos \in \C{U}$
that
\begin{equation}   \label{eq:def-simpl-ex}
   \vek{R} \vek{u} = 
     [{\phi(\mu_1)}, {\phi(\mu_2)}, {\phi(\mu_3)}]^\trpos = \vek{\phi} \in \D{R}^3
     = \C{F} = \D{R}^{\C{M}},
\end{equation}
where $\phi(\mu_i) = \sum_{j=1}^4 r_j(\mu_i) u_j$ --- a \emph{weighted} average
of $r(\cdot,\mu_i)$ --- is a scalar function $\phi\in\C{F} = \D{R}^{\C{M}}$
in the linear space $\C{F}$ of scalar functions $\D{R}^{\C{M}}$ on the parameter set $\C{M}$.
If we denote the function of \refeq{eq:def-r-map} in this case by $\vek{r}(\cdot)$,
which for every $\mu\in\C{M}$ is an element $\vek{w}:=\vek{r}(\mu)\in\C{U} =\D{R}^4$,
then the weighted average $\vek{\phi}\in\C{R} = \D{R}^3$ in \refeq{eq:def-simpl-ex}
obviously satisfies $\phi(\mu_i) = \vek{r}(\mu_i)^\trpos \vek{u}$, so that
\begin{equation}   \label{eq:def-simpl-ex-2}
   \vek{R} \vek{u} = [\vek{r}(\mu_1)^\trpos \vek{u}, \vek{r}(\mu_2)^\trpos \vek{u},
      \vek{r}(\mu_3)^\trpos \vek{u}]^\trpos = \ip{\vek{r}(\cdot)}{\vek{u}}_{\C{U}} =
     [{\phi(\mu_1)}, {\phi(\mu_2)}, {\phi(\mu_3)}]^\trpos = \vek{\phi} .
\end{equation}
Obviously, knowing $\vek{R}$ is the same as knowing $\vek{R} \vek{u}$ for every
$\vek{u}\in\C{U} =\D{R}^4$ --- actually a basis in $\C{U}$ would suffice --- which
in turn is the same as knowing $\ip{\vek{r}(\cdot)}{\vek{u}}_{\C{U}}$ for every
$\vek{u}\in\C{U}$.

The point to take away from this simple example is that the parametric function
$r(x,\mu)$, where for each parameter value $\mu\in\C{M}$ one has
$r(\cdot,\mu)\in\C{U}$ in some linear space $\C{U}$ --- of functions on $\C{X}$
in this case --- is equivalent to a linear map 
\[ \vek{R}:\C{U}\to\C{F} \]
into a space $\C{F}\subseteq\D{R}^{\C{M}}$ of scalar functions on the parameter set $\C{M}$.

It is now easy to see how to generalise this further to cases where the set
$\C{X}$ or $\C{M}$ or both have infinitely many values, and even further to a
case where the vector space of functions $\C{U}$ just has an inner product, say
given by some integral, so that for $u, v \in \C{U}$ one has
\[ \ip{u}{v}_{\C{U}} = \int_{\C{X}} u(x) v(x) \,\ops{m}(\di x) \]
with some measure $\ops{m}$ on $\C{X}$.  Then for each parameter $\mu\in\C{M}$ one
has $r(\cdot,\mu)\in\C{U}$, a function on $\C{X}$, or in other words
an element of the linear space $\C{U}$.  In this case one defines the linear map
\[   \tilde{R}:\C{U}\to\C{F}\subseteq\D{R}^{\C{M}} \] as  
\[
  \tilde{R}: u \mapsto \int_{\C{X}} u(x) r(x,\mu) \,\ops{m}(\di x) = 
     \ip{r(\cdot,\mu}{u}_{\C{U}} =: \phi(\mu) \in\C{F} \subseteq\D{R}^{\C{M}},
\]
which is a linear map from $\C{U}$ onto a linear space of scalar functions
$\phi\in\C{F}\subseteq\D{R}^{\C{M}}$ on the parameter set $\C{M}$.

This then is almost the general situation, where one views $r:\C{M}\to\C{V}$
as a map from the parameters $\mu\in\C{M}$, where $\C{M}$ may be some arbitrary
set, into a topological vector space $\C{V}$. One then defines a linear map
\[ \tilde{R}:\C{V}^* \to \C{F} \subseteq\D{R}^{\C{M}} \]
from the dual space $\C{V}^*$ onto a space of scalar functions $\C{F}$ on $\C{M}$ by
\[
    \C{V}^* \ni u \mapsto \tilde{R}u = \bkt{r(\mu)}{u}_{(\C{V},\C{V}^*)} =: \phi(\mu) \in 
      \C{F} \subseteq\D{R}^{\C{M}},
\]
where $\bkt{\cdot}{\cdot}_{(\C{V},\C{V}^*)}$ is the duality pairing between
$\C{V}$ and its dual space $\C{V}^*$.  For the following exposition of the main
ideas we shall take a slightly less general situation by assuming for the
sake of simplicity that
the linear space $\C{V}$ is in fact a separable Hilbert space with an inner product
$\ip{\cdot}{\cdot}_{\C{U}}$, and use this in the usual manner to identify it with its dual.

\paragraph{Associated linear map:}
So with a vector-valued map $r:\C{M}\to\C{V}$, one defines the corresponding
associated linear map $\tilde{R}:\C{V}\to\C{F}$ as
\begin{equation}     \label{eq:V-1}
\forall v \in \C{V}: \tilde{R} v = \ip{r(\mu)}{v}_{\C{U}} =: \phi(\mu) \in \C{F}
                  \subseteq\D{R}^{\C{M}}.
\end{equation}
Obviously only the Hilbert subspace $\C{U}=\cl (\spn r(\C{M}))\subseteq\C{V}$
actually reached by the map $r$ is interesting, whereas
$\C{U}^\perp = \ker \tilde{R}\subseteq\C{V}$ is not.  
Hence from now on we shall only look at $r:\C{M}\to\C{U}$, and additionally
assume that $\C{U} = \cl (\spn r(\C{M}))$, or in other words, that the vectors
$\{ r(\mu) \mid \mu \in \C{M} \} = r(\C{M})$ form a total set in $\C{U}$.
The map $\tilde{R}$ is thus formally redefined as
\begin{equation}     \label{eq:U-1}
  \tilde{R}: \C{U} \to \D{R}^{\C{M}}.
\end{equation}
Again, in the linear space $\D{R}^{\C{M}}$ of \emph{all} scalar functions on $\C{M}$,
only the part $\C{F} = \im \tilde{R} = \tilde{R}(\C{U})$ is interesting.

Allow here a little digression, to point out similarities and
analogies to other connected concepts.  First, on the parameter set $\C{M}$,
where up to now no additional mathematical structure was used,
we now have the linear space $\C{F}$.  This can be viewed as a first step
to introduce some kind of ``co-ordinates'' on the set $\C{M}$, and is in line
with many other constructs where potentially complicated sets are characterised
by algebraic constructs, such as groups or vector spaces  for e.g.\ homology
or cohomology.  Even if from the outset the parameter set $\C{M}\subseteq \D{R}^m$
is given as some subset of some $\D{R}^m$ and therefore has already coordinates,
these may not be \emph{good} ones, and as we shall see, it may be worthwhile
to contemplate \emph{re-parametrisations}, i.e.\ choosing some $\phi_k \in \C{F}$
as ``co-ordinates''.  These real valued functions are in general of course not
``real co-ordinates'', as they only distinguish what is being felt by the parametric
object $r$.

\paragraph{Reproducing kernel Hilbert space:}  The second concept to touch on
comes from the idea to use the function space $\C{F}$ in place of $\spn r(\C{M})$:
As is easy to see, the map in \refeq{eq:U-1} is injective, hence invertible on its
image $\C{F} = \im \tilde{R} = \tilde{R}(\C{U})$, and this may be used
to define an inner product on $\C{F}$ as
\begin{equation}     \label{eq:V}
\forall \phi, \psi \in \C{F} \quad \ip{\phi}{\psi}_{\C{R}} :=
     \ip{\tilde{R}^{-1} \phi}{\tilde{R}^{-1} \psi}_{\C{U}},
\end{equation}
and to denote the completion of $\C{F}$ with this inner product by $\C{R} = \cl \C{F}
\subseteq \D{R}^{\C{M}}$.
One immediately obtains that $\tilde{R}^{-1}$ is a bijective isometry between $\spn \im r$
and $\C{F}$, hence extends to a \emph{unitary} map $\bar{R}^{-1}$ between $\C{U}$ and $\C{R}$,
and the same hold for $\tilde{R}$, the extension being denoted by $\bar{R}$.

Given the maps $r:\C{M}\to\C{U}$ and $\bar{R}:\C{U}\to\C{R}$, one may define the
\emph{reproducing kernel} \citep{berlinet, Janson1997}
given by $\vkappa(\mu_1, \mu_2) := \ip{r(\mu_1)}{r(\mu_2)}_{\C{U}}$.
It is straightforward to verify that $\vkappa(\mu,\cdot)\in\C{F}\subseteq\C{R}$,
and $\spn \{ \vkappa(\mu,\cdot)\;\mid\; \mu\in\C{M} \}=\C{F}$,
as well as the reproducing property
$\phi(\mu) = \ip{\vkappa(\mu,\cdot)}{\phi}_{\C{R}}$ for all $\phi\in\C{F}$.
Another way of stating this reproducing property is to say that the linear map 
$\C{R}\ni\phi(\mu_1) \mapsto \ip{\vkappa(\mu_2,\cdot)}{\phi}_{\C{R}}
 = \phi(\mu_2) \in\C{R}$ for all $\phi\in\C{R}$
is the identity $I_{\C{R}}$  on $\C{R}$.
An abstract way of putting this using the adjoint $\bar{R}^* = \bar{R}^{-1}$ of the
unitary map $\bar{R}$ would be to note that that map is in fact 
$\bar{R} \bar{R}^* = \bar{R} \bar{R}^{-1} = I_{\C{R}}$.

With the \emph{reproducing kernel Hilbert space} (RKHS) $\C{R}$
one can build a first representation and thus obtain a relevant
``co-ordinate system'' for $\C{M}$.  As $\C{U}$ is separable, it has a Hilbert basis
or complete orthonormal system (CONS) $\{y_k\}_{k\in\D{N}}$.
As $\bar{R}$ is unitary, the set $\{ \vphi_k = \bar{R} y_k \}_{k\in\D{N}}$
is a CONS in $\C{R}$.

With this, the unitary operator $\bar{R}$, 
its adjoint or inverse $\bar{R}^*=\bar{R}^{-1}$,
and the parametric element $r(\mu)$ become \citep{hgmRO-1-2018}
\begin{align} \label{eq:VII0}
  \bar{R} &= \sum_k \vphi_k \otimes y_k;  \quad \text{i.e. } \quad 
  \bar{R}(u)(\cdot) = \sum_k \bkt{y_k}{u}_{\C{U}} \vphi_k(\cdot), \\  \label{eq:VII0-1}
  \bar{R}^* =\bar{R}^{-1} &= \sum_k y_k \otimes \vphi_k;\quad
   r(\mu) = \sum_k \vphi_k(\mu) y_k = \sum_k \vphi_k(\mu)\, \bar{R}^* \vphi_k .
\end{align}
Observe that the relations \refeq{eq:VII0} and \refeq{eq:VII0-1} exhibit 
the tensorial nature of the representation mapping.
One sees that \emph{model reductions} may be achieved by choosing only subspaces of $\C{R}$,
i.e.\ spanned by a---typically finite---subset of the CONS $\{\vphi_k\}_{k}$.
Furthermore, the representation of $r(\mu)$ in \refeq{eq:VII0-1} is \emph{linear}
in the new ``parameters'' $\vphi_k$.

\paragraph{Coherent states:}  The third concept one should mention in this context
is the one of \emph{coherent states}, e.g. see \citep{AliAntGazeau2014, AntBagGazeau2018}.
In this development from quantum theory, these quantum states were initially introduced
as eigenstates of certain operators, and the name refers originally to their high coherence,
minimum uncertainty, and quasi classical behaviour.  What is important here is that
the idea has been abstracted, and represents overcomplete sets of vectors or frames
$\{ r(\mu) \mid \mu\in\C{M} \}$
in a Hilbert space $\C{U}$, which depend on a parameter $\mu\in\C{M}$ from a locally
compact measure space.  This space often has more structure, e.g.\ a Lie group,
and the coherent states are connected with group representations in the unitary group
of $\C{U}$, i.e.\ if $\mu \mapsto U(\mu) \in \E{L}(\C{U})$ is a unitary representation,
the coherent states may be defined by $r(\mu) = U(\mu) w$ for some $w\in\C{U}$.
There are usually further requirements like weak continuity for
the map $\C{M}\ni\mu \mapsto r(\mu)\in\C{U}$, and that these coherent states form
a \emph{resolution of the identity}, in that one has (weakly)
\[
  I_{\C{U}} = \int_{\C{M}} r(\mu)\otimes r(\mu) \, \vpi(\di \mu) ,
\] 
where $\vpi$ is a measure on $\C{M}$---naturally defined on some $\sigma$-algebra
of subsets of $\C{M}$, a detail which needs no further mention here.  
We shall leave this topic here, and
come back to similar representations later, but note in passing the tensor
product structure under the integral.  The above requirement of
the resolution of the identity may sometimes be too
strong, and one often falls to the case of RKHS discussed above.

\paragraph{Reduced models:} Assume now that $\C{M}\ni\mu \mapsto r_a(\mu)\in\C{U}$
is an approximate or \emph{reduced order model} (ROM) of $r(\mu)$.  One possibility
of producing such a ROM was already mentioned above by letting the sum in
\refeq{eq:VII0-1} run over fewer terms.  The ROM $r_a(\mu)$ thus has an
associated linear map $\bar{R}_a$.  As the associated linear maps carry all the
relevant information, the analysis of both the original parametric object $r(\mu)$,
and the comparison and analysis of accuracy of the approximation $r_a(\mu)$ can
be carried out in terms of the associated linear maps $\bar{R}$ and $\bar{R}_a$.
In the present setting $\bar{R}$ is unitary, so $\bar{R}_a$ can be judged by how
well it approximates that unitary mapping.  In the next \refS{correlat}, where 
a second inner product will be introduced on the space $\D{R}^{\C{M}}$ of
scalar functions of $\C{M}$, this will be even more pronounced, as it will offer
the possibility of deciding which CONS or other complete sets in subspaces
of $\D{R}^{\C{M}}$ are advantageous for ROMs.


%
%
%
%


%

\section{Correlation and Representation}  \label{S:correlat}
In what was detailed up to now in the previous \refS{parametric} with regard
to the RKHS, was that the structure of the Hilbert space was carried reproduced
on the subspace $\C{R}\subseteq \D{R}^{\C{M}}$ of the full function space.
In the remarks about coherent states one could already see an additional structure,
namely a measure $\vpi$ on $\C{M}$.
This measure structure can be used to define the subspace $\C{A}:=\Lp_0(\C{M},\vpi)$ of
measurable functions, as well as its Hilbert subspace of square-integrable functions
$\C{Q}:=\Lp_2(\C{M},\vpi)$ with associated inner product
\[
  \ip{\phi}{\psi}_{\C{Q}} := \int_{\C{M}} \phi(\mu) \psi(\mu) \; \vpi(\di \mu).
\]
We shall simply assume here that there is a Hilbert space $\C{Q}\subseteq \D{R}^{\C{M}}$
of functions with inner product $\ip{\cdot}{\cdot}_{\C{Q}}$, which may or may not come
from an underlying measure space.  The associated linear map $\tilde{R}:\C{U}\to\C{R}$,
essentially defined in \refeq{eq:V-1} with range the RKHS $\C{R}$, will now be seen
as a map  $R:\C{U}\to\C{Q}$ into the Hilbert space $\C{Q}$, i.e.\ with a different
range with different inner product $\ip{\cdot}{\cdot}_{\C{Q}}$ from 
the RKHS inner product $\ip{\cdot}{\cdot}_{\C{R}}$ on $\C{R}$.  One may view
this inner product as a way to tell what is important in the parameter set
$\C{M}$: functions $\phi$ with large $\C{Q}$-norm are considered more important than
those where this norm is small.  The map $R:\C{U}\to\C{Q}$ is thus generally not
unitary any more, but for the sake of simplicity, we shall assume that it is
a densely defined closed operator, see e.g.~\citep{DautrayLions3}.  As it may be
only densely defined, it is sometimes a good idea to define $R$ through a densely
defined bilinear form in $\C{U}\otimes\C{Q}$:
\begin{equation}   \label{eq:IX-a}
   \forall u\in \dom R, \phi \in \C{Q}: \ip{Ru}{\phi}_{\C{Q}} :=
   \ip{\ip{r(\cdot)}{u}_{\C{U}}}{\phi}_{\C{Q}}.
\end{equation}

Following \citep{kreeSoize86, hgmRO-1-2018, hgmRO-2-2018, hgm-3-2018}, one now obtains
a densely defined map $C$ in $\C{U}$ through the densely defined bilinear form,
in line with \refeq{eq:IX-a}:
\begin{equation}   \label{eq:IX}
   \forall u, v:\quad \ip{Cu}{v}_{\C{U}} := \ip{Ru}{Rv}_{\C{Q}} .
\end{equation}
The map $C=R^* R$ --- observe that now the adjoint is w.r.t.\ the $\C{Q}$-inner
product --- may be called the \emph{``correlation''} operator, and is by
construction self-adjoint and positive, and if $R$ is bounded resp.\ continuous, so is $C$.

In the above case that the $\C{Q}$-inner product comes from a measure, 
one has from \refeq{eq:IX}
\[  
  \ip{Cu}{v}_{\C{U}} = 
  \int_{\C{M}} \ip{r(\mu)}{u}_{\C{U}}\ip{r(\mu)}{v}_{\C{U}} \; \vpi(\di \mu),\text{ i.e. }
  C = R^* R = \int_{\C{M}} r(\mu) \otimes r(\mu) \; \vpi(\di \mu).
\]
This is reminiscent of what was required for coherent states.  But it also shows
that if $\vpi$ were a probability measure --- i.e.\ $\vpi(\C{M})=1$ --- with the
usual expectation operator 
\[ \EXP{\phi} := \int_{\C{M}} \phi(\mu) \; \vpi(\di \mu), \]
then the above would be really the familiar correlation operator
\citep{kreeSoize86, hgm-3-2018}
$\EXP{r\otimes r}$ of the $\C{U}$-valued \emph{random variable} (RV) $r$, therefore
from now on we shall simply refer to $C$ as the correlation operator, even
in the general case not based on a probability measure.

The fact that the correlation operator is self-adjoint and positive
implies that its spectrum $\sigma(C)\subseteq \D{R}^+$ is real and non-negative.
This will be used when analysing it with any of the versions of the spectral
theorem for self-adjoint operators (e.g.~\citep{DautrayLions3}).  
The easiest and best known version of this is for finite dimensional maps.

\paragraph{Finite dimensional beginnings:}
So let us return to the simple example at the beginning of \refS{parametric}
where the associated linear map can be represented by a matrix $\vek{R}$.
If we remember the each row $\vek{r}^\trpos(\mu_j)$ is the value for the
vector $\vek{r}(\mu)$ for one particular
$\mu\in\C{M}$, we see that the matrix can be written as
\[
   \vek{R} = [\vek{r}(\mu_1),\dots,\vek{r}(\mu_j),\dots ]^\trpos ,
\]
and that the rows are just ``snapshots'' for different values $\mu_j$.
What is commonly done now is the so-called method
of proper orthogonal decomposition (POD) to produce a ROM.

The matrix $\vek{R}$ --- to generalise a bit, assume it of size $m\times n$ ---
can be decomposed according to its \emph{singular value decomposition} (SVD)
\begin{equation}  \label{eq:SVD-R-fdim}
   \vek{R} = \vek{\Phi}\vek{\Sigma}\vek{V}^\trpos = 
   \sum_{k=1}^{\min(m,n)} \vsigma_k\, \vek{\phi}_k \otimes \vek{v}_k,
\end{equation}
where the matrices $\vek{\Phi}=[\vek{\phi}_k]$ and $\vek{V}=[\vek{v}_k]$ are
orthogonal with unit length orthogonal columns --- right and left singular vectors ---
$\vek{\phi}_k$ resp.\ $\vek{v}_k$, and $\vek{\Sigma} = \diag(\vsigma_k)$ is
diagonal with non-negative diagonal elements $\vsigma_k$, the singular values.
For clarity, we arrange the singular values in a decreasing sequence,
$\vsigma_1 \ge \vsigma_2 \ge \dots \ge 0$.
It is well known that this decomposition is connected with the eigenvalue
or spectral decomposition of the correlation
\begin{equation}  \label{eq:spec-C-fdim}
  \vek{C} = \vek{R}^\trpos \vek{R} = \vek{V}\vek{\Sigma}\vek{\Phi}^\trpos
  \vek{\Phi}\vek{\Sigma}\vek{V}^\trpos = \vek{V}\vek{\Sigma}^2\vek{V}^\trpos =
  \sum_{k=1}^{\min(m,n)}  \vsigma_k^2 \, \vek{v}_k\otimes \vek{v}_k,
\end{equation}
with eigenvalues $\vsigma_k^2$, eigenvectors $\vek{v}_k$, and its companion
\begin{equation}  \label{eq:spec-CQ-fdim}
  \vek{C}_{\C{Q}} := \vek{R} \vek{R}^\trpos  = \vek{\Phi}\vek{\Sigma}^2\vek{\Phi}^\trpos =
  \sum_{k=1}^{\min(m,n)}  \vsigma_k^2 \, \vek{\phi}_k\otimes \vek{\phi}_k,
\end{equation}
with the same eigenvalues, but eigenvectors $\vek{\phi}_k$.
The representation is based on $\vek{R}^\trpos$, and its accompanying POD or
\emph{Karhunen-Loève} decomposition:
\begin{equation}  \label{eq:KLE-fdim}
  \vek{R}^\trpos = \sum_{k=1}^{\min(m,n)} \vsigma_k\, \vek{v}_k \otimes \vek{\phi}_k ,\qquad
  \vek{r}(\mu_j) = \sum_{k=1}^{\min(m,n)} \vsigma_k\, \vek{v}_k \otimes \vek{\phi}_k(\mu_j),
\end{equation}
where $\vek{\phi}_k(\mu_j)=\phi_k^j$, and 
$\vek{\phi}_k = [\phi_k^1,\dots,\phi_k^j,\dots]^\trpos $. 

The second expression in \refeq{eq:KLE-fdim} is a representation for $\vek{r}(\mu)$,
and that is the purpose of the whole exercise.  Similar expressions may be used as
approximations.  It clearly exhibits the tensorial
nature of the representation, which is also evident in the expressions \refeq{eq:SVD-R-fdim},
\refeq{eq:spec-C-fdim}, and \refeq{eq:spec-CQ-fdim}.  One sees here that this is just the
$j$-th column of $\vek{R}^\trpos$, so that with the canonical basis in $\C{Q}=\D{R}^m$,
$\vek{e}_j^{(m)} = [\updelta_{ij}]^\trpos$ with the Kronecker-$\updelta$, that expression
becomes just
\begin{equation}  \label{eq:KLE-fdim-appr}
  \vek{r}(\mu_j) = \vek{R}^\trpos \vek{e}_j^{(m)}; \quad \text{ and } \quad
  \vek{r}\approx\vek{r}_a = \vek{R}^\trpos \vek{\psi}
\end{equation}
by taking other vectors $\vek{\psi}$ in $\C{Q}=\D{R}^m$ to give weighted
averages or interpolations.

The general picture which emerges is that the matrix $\vek{R}$ is a kind 
of ``square root''  --- or more precisely factorisation --- of
the correlation $\vek{C}=\vek{R}^\trpos\vek{R}$, and that the left part of
this factorisation is used for reconstruction resp.\ representation.
In any other factorisation like
\begin{equation}  \label{eq:fact-fdim}
   \vek{C} = \vek{B}^\trpos\vek{B}, \quad \text{ with } \quad 
   \vek{B}:\C{U}\to\C{H},
\end{equation}
where $\vek{B}$ maps into some other space $\C{H}$; the map $\vek{B}$ will
necessarily have essentially the same singular values $\vsigma_k$ and right singular
vectors $\vek{v}_k$ as $\vek{R}$, and can now be used to have a representation
or reconstruction of $\vek{r}$ on $\C{H}$ via
\begin{equation}  \label{eq:fact-fdim-h-appr}
   \vek{r} \approx \vek{B}^\trpos \vek{h} \quad \text{ for some } \quad \vek{h}\in\C{H}.
\end{equation}
A popular choice is to use the Choleski-factorisation $\vek{C}=\vek{L}\vek{L}^\trpos$
of the correlation into two triangular matrices, and then take $\vek{B}^\trpos=\vek{L}$
for the reconstruction.

As we have introduced the correlation's spectral factorisation in \refeq{eq:spec-C-fdim},
some other factorisations come to mind, although they may be mostly of
theoretical value:
\begin{equation}  \label{eq:fact2-fdim}
   \vek{C} = \vek{B}^\trpos\vek{B} = (\vek{V}\vek{\Sigma})(\vek{V}\vek{\Sigma})^\trpos
   = (\vek{V}\vek{\Sigma}\vek{V}^\trpos)(\vek{V}\vek{\Sigma}\vek{V}^\trpos)^\trpos,
\end{equation}
where then the reconstruction map is $\vek{B}^\trpos = (\vek{V}\vek{\Sigma})$ or
$\vek{B}^\trpos = (\vek{V}\vek{\Sigma}\vek{V}^\trpos)$.  Obviously, in the second case
the reconstruction map is symmetric $\vek{B}^\trpos=\vek{B}=\vek{C}^{1/2}$, and is
actually the true square root of the correlation $\vek{C}$.

Other factorisation can come from looking at the companion $\vek{C}_{\C{Q}}$ in
\refeq{eq:spec-CQ-fdim}.  Any factorisation $\vek{F}:\C{Z}\to\C{Q}$
or approximate factorisation $\vek{F}_a$ of
\begin{equation}  \label{eq:fact-Q-fdim}
   \vek{C}_{\C{Q}} = \vek{F}\vek{F}^\trpos \approx \vek{F}_a\vek{F}_a^\trpos
\end{equation}
is naturally a factorisation or approximate factorisation of the correlation
\begin{equation}  \label{eq:fact-CQ-fdim}
   \vek{C} = \vek{W}^\trpos \vek{W}
    \approx \vek{W}_a^\trpos\vek{W}_a, \quad \text{ with } \quad
   \vek{W} = \vek{F}^\trpos\vek{\Phi}\vek{V}^\trpos \text{ and } 
   \vek{W}_a = \vek{F}_a^\trpos \vek{\Phi}\vek{V}^\trpos ,
\end{equation}
where $\vek{V}$ and $\vek{\Phi}$ are the left and right singular vectors
--- see \refeq{eq:SVD-R-fdim} --- of the associated map
$\vek{R}$ resp.\  the eigenvectors of the correlation $\vek{C}$ in 
\refeq{eq:spec-C-fdim} and its companion $\vek{C}_{\C{Q}}$ in \refeq{eq:spec-CQ-fdim}.
A new ROM representation can now be found for $\vek{z}\in\C{Z}$ via
\begin{equation}  \label{eq:ROM-Q-fdim}
\vek{r} \approx \vek{r}_a = \vek{W}_a^\trpos \vek{z} =  
\vek{V}\vek{\Phi}^\trpos \vek{F}_a \vek{z}.
\end{equation}
One last observation here is important: the expressions for $\vek{r}$ resp.\ one
of its ROMs $\vek{r}_a$ are \emph{linear} in the newly introduced parameters
or ``co-ordinates'' $\vek{\phi}_k$ in \refeq{eq:KLE-fdim}, resp.\ $\vek{\psi}$ 
in \refeq{eq:KLE-fdim-appr}, resp.\ $\vek{h}$ in \refeq{eq:fact-fdim-h-appr}
and \refeq{eq:ROM-fdim},  as well as $\vek{z}$ in \refeq{eq:ROM-Q-fdim};
which is an important requirement in many numerical methods.

\paragraph{Reduced order models --- ROMs:}
As has become clear now, and was mentioned before, 
that approximations or ROMs $\vek{r}_a(\mu)$ to the full
model $\vek{r}(\mu) \approx \vek{r}_a(\mu)$ produce associated maps $\vek{R}_a$,
which are approximate factorisations of the correlation:
\[ \vek{C} \approx \vek{R}_a^\trpos \vek{R}_a . \]
This introduces different ways of judging how good an approximation is.
If one looks at the difference between the full model $\vek{r}(\mu)$
and ist approximation $\vek{r}_a(\mu)$ as a residual, and computes weighted
versions of it
\begin{equation}   \label{eq:r-ra-diff}
   \ip{\vek{r}(\cdot) - \vek{r}_a(\cdot)}{u}_{\C{U}} = (\vek{R}-\vek{R}_a) \vek{u}
   = \vek{R} \vek{u} - \vek{R}_a  \vek{u} ,
\end{equation}
then this is just the difference linear map $\vek{R}-\vek{R}_a$ applied
to the weighting vector $\vek{u}$.  In \refeq{eq:KLE-fdim} is was shown that
$\vek{r}(\cdot) = \sum_{k=1}^{\min(m,n)} \vsigma_k\, \vek{v}_k \otimes \vek{\phi}_k(\cdot)$
is a representation.  As usual, one may now approximate such an expressions
by leaving out terms with small or vanishing singular values, say using only 
$\vsigma_1, \dots,\vsigma_\ell$, getting an approximation of rank $\ell$ ---
this also means that the associated linear map $\vek{R}_a$ in \refeq{eq:KLE-fdim}
has rank $\ell$.  As is well known \citep{Hackbusch_tensor}, this is the best $\ell$-term
approximation in the norms of $\C{U}$ and $\C{Q}$.
But from \refeq{eq:r-ra-diff} one may gather that the error can also be
described through the difference $\vek{R}-\vek{R}_a$.
As error measure one may take the norm of that difference, and, depending
on which norm one chooses, the error is then in  --- this example approximation
 --- $\vsigma_{\ell+1}$ in the operator norm, or
$\sum_{k=\ell+1}^{\min(m,n)}\vsigma_k$ in the trace- resp.\ nuclear norm, or
$\sqrt{\sum_{k=\ell+1}^{\min(m,n)}\vsigma_k^2}$ in the Frobenius- resp. Hilbert-Schmidt norm.

On the other hand, different approximations or ROMs can now be obtained by starting with
an approximate factorisation
\begin{equation}  \label{eq:fact-a-fdim}
   \vek{C} \approx \vek{B}_a^\trpos\vek{B}_a,
\end{equation}
and introducing a ROM via 
\begin{equation}  \label{eq:ROM-fdim}
\vek{r} \approx \vek{r}_a = \vek{B}_a^\trpos \vek{h} .
\end{equation}
Such a representing linear map $\vek{B}$, may, e.g.\ via its SVD,
be written as a sum of tensor products, and approximations $\vek{B}_a$ are often
lower rank expressions, directly reflected in a reduced sum for the
tensor products.  As will become clearer at the end of this section,
the bilinear forms \refeq{eq:IX-a} resp.\ \refeq{eq:IX} can sometimes split
into multi-linear forms, thus enabling the further approximation of
$\vek{B}_a$ through \emph{hierarchical tensor products} \citep{Hackbusch_tensor}.

\paragraph{Infinite dimensional continuation --- discrete spectrum:}
For the cases where both $\C{U}$ and $\C{Q}$ are infinite dimensional, the
operators $R$ and $C$ live on infinite dimensional spaces, and the spectral
theory gets a bit more complicated.  We shall distinguish some simple
cases.  After finite dimensional resp.\ finite rank operators just treated in
matrix form, the next simplest
case is certainly the case when the associated linear map $R$ and the
correlation operator $C=R^* R$ has a discrete spectrum, e.g.\ if $C$ is compact,
or a function of a compact operator, like for example its inverse.
In this case the spectrum is discrete (e.g.~\citep{DautrayLions3}),
and in the case of a compact operator the non-negative eigenvalues $\lambda_k$
of $C$ may be arranged as a decreasing sequence 
$\infty>\lambda_1\ge\lambda_2\ge\dots\ge 0$ with only possible accumulation point
the origin.  It is not uncommon when dealing with random fields that
$C$ is a nuclear or trace-class operator, i.e.\ an operator which satisfies
the stronger requirement $\sum_k \lambda_k < \infty$.  
The spectral theorem for a an operator with purely discrete spectrum takes the form
\begin{equation}  \label{eq:XIII}
   C = R^* R = \sum_{k=1}^\infty \lambda_k \, (v_k \otimes v_k) ,
\end{equation}
where the eigenvectors $\{v_k\}_k \subset \C{U}$ form a CONS in $\C{U}$.
Defining a new corresponding CONS $\{s_k\}_k$ in $\C{Q}$ via
 $\lambda_k^{1/2} s_k := R v_k$, one obtains the
singular value decomposition of $R$ and $R^*$
with singular values $\vsigma_k=\lambda_k^{1/2}$:
\begin{multline}   \label{eq:XIV}
R = \sum_{k=1}^\infty \vsigma_k (s_k \otimes v_k)\,;   \quad \text{i.e. } \quad 
  R(u)(\cdot) = \sum_{k=1}^\infty \vsigma_k \ip{v_k}{u}_{\C{U}} s_k(\cdot), \quad
R^* = \sum_{k=1}^\infty \vsigma_k (v_k \otimes s_k)\,; \\
r(\mu) =  \sum_{k=1}^\infty \vsigma_k \, s_k(\mu) v_k =
 \sum_{k=1}^\infty s_k(\mu)\, R^* s_k,
\text{ as } \; R^*s_k = \vsigma_k \, v_k.
\end{multline}
It is not necessary to repeat in this setting of compact maps all the different
factorisations considered in the preceding paragraphs, and especially their approximations,
which will be usually finite dimensional as they are made to be used for actual
computations, e.g.\  the approximations will usually involve only finite portions
of the infinite series in \refeq{eq:XIII} and \refeq{eq:XIV}, which means that
the induced linear maps have finite rank and essentially become finite dimensional,
so that the preceding paragraphs apply practically verbatim.

But one consideration is worth to follow up further.
In infinite dimensional Hilbert spaces, self-adjoint operators may have
a continuous spectrum, e.g.~\citep{DautrayLions3};  this is what is usually
the case when homogeneous random fields or stationary stochastic processes
have to be represented,  This means that the expressions developed for purely
discrete spectra in \refeq{eq:XIII} and \refeq{eq:XIV} are not general enough.
These expressions are really generalisations of the last equalities in
\refeq{eq:spec-C-fdim} and \refeq{eq:SVD-R-fdim}; but is is possible 
to give meaning to the matrix equalities in those equations, which simultaneously
cover the case of a continuous spectrum.

\paragraph{In infinite dimensions --- non-discrete spectrum:}
To this end we introduce the so called multiplication operator: Let
$\Lp_2(\C{T})$ be the usual Hilbert space on some locally compact
measure space $\C{T}$, and let $\gamma\in\Lp_\infty(\C{T})$ be an essentially
bounded function.  Then the map 
\[ M_{\gamma}:\Lp_2(\C{T})\to\Lp_2(\C{T}); \qquad M_{\gamma}:\xi(t) \to \gamma(t)\xi(t) \]
for $\xi \in \Lp_2(\C{T})$ is a bounded operator $M_{\gamma}\in\E{L}(\C{X})$
on $\Lp_2(\C{T})$.  Such a multiplication operator is the direct analogue of a
diagonal matrix in finite dimensions.

%

Using such a multiplication operator, one may introduce
a formulation of the spectral decomposition different from \refeq{eq:XIII}
which does not require $C$ to be compact \citep{DautrayLions3},
$C$ resp.\ $R$ do not even have to be continuous resp.\ bounded:
\begin{equation}  \label{eq:ev-mult}
  C = R^* R = V M_{\gamma} V^*,
\end{equation}
where $V:\Lp_2(\C{T})\to\C{U}$ is unitary between some $\Lp_2(\C{T})$ 
on a measure space $\C{T}$ and $\C{U}$.  In case $C$ is continuous resp.\ bounded,
one has $\gamma\in\Lp_\infty(\C{T})$.  As $C$ is positive, the function
$\gamma$ is non-negative ($\gamma(t) \ge 0$ a.e.\ for $t\in\C{T}$).
This covers the previous case
of operators with purely discrete spectrum if the function $\gamma$
is a step function and
takes only a discrete (countable) set of values --- the eigenvalues.
This theorem is actually quite well known in the special case that
$C$ is the correlation operator of a stationary stochastic process
--- an integral operator where the kernel is the correlation function;
in this case $V$ is the \emph{Fourier transform}, and $\gamma$ is
known as the \emph{power spectrum}.

\paragraph{General factorisations:}
To investigate the analogues of further factorisations of $R$,
$C=R^* R$, and its companion $C_{\C{Q}} = R R^*$, we need the
SVD of $R$ and $R^*$.  They derive generally in the same manner as
for the finite dimensional case from the spectral factorisations of $C$ in
\refeq{eq:ev-mult} and a corresponding one for its companion
\begin{equation}  \label{eq:ev-cq-mult}
  C_{\C{Q}} = R R^* = \Phi M_{\gamma} \Phi^* 
\end{equation}
with a unitary $\Phi:\Lp_2(\C{T}_*)\to\C{Q}$ between some $\Lp_2(\C{T}_*)$ 
on a measure space $\C{T}_*$ and $\C{Q}$.  Here in \refeq{eq:ev-cq-mult},
and in \refeq{eq:ev-mult}, the multiplication operator $M_{\gamma}$ plays
the role of the diagonal matrix $\vek{\Sigma}^2$ in \refeq{eq:spec-C-fdim}
and \refeq{eq:spec-CQ-fdim}. For the SVD of $R$ one needs its square root,
and as $\gamma$ is non-negative, this is simply given by
$M_{\gamma}^{1/2} = M_{\sqrt{\gamma}}$, i.e.\ multiplication by $\sqrt{\gamma}$.
Hence the  SVD of $R$ and $R^*$ is given by
\begin{equation}  \label{eq:ev-mult-svd}
R = \Phi M_{\sqrt{\gamma}} V^*,\quad R^* = V M_{\sqrt{\gamma}} \Phi^*.
\end{equation}

These are all examples of a general factorisation
$C = B^* B$, where $B:\C{U}\to\C{H}$ is a map to a Hilbert space $\C{H}$
with all the properties demanded from $R$---see the beginning of this section.
It can be shown \citep{hgmRO-1-2018} that
any two such factorisations $B_1:\C{U}\to\C{H}_1$ and $B_2:\C{U}\to\C{H}_2$ 
with $C=B_1^*B_1=B_2^*B_2$ are
\emph{unitarily equivalent} in that there is a unitary map $X_{21}:\C{H}_1\to\C{H}_2$
such that $B_2 = X_{21} B_1$.  Equivalently, each such factorisation is unitarily
equivalent to $R$, i.e.\  there is a unitary $X:\C{H}\to\C{Q}$ such that
$R= X B$.

Analogues of the factorisations considered in \refeq{eq:fact2-fdim} are
\begin{equation}  \label{eq:fact2}
   C = B^* B = (V M_{\sqrt{\gamma}})(V M_{\sqrt{\gamma}})^*
   = (V M_{\sqrt{\gamma}}V^*)(V M_{\sqrt{\gamma}V}^*)^*,
\end{equation}
where again $C^{1/2} = V M_{\sqrt{\gamma}}V^*$ is the square root of $C$.

And just as in the case of the factorisations of $\vek{C}_{\C{Q}}$ considered
in \refeq{eq:fact-Q-fdim} and the resulting factorisation of $\vek{C}$ in
\refeq{eq:fact-CQ-fdim}, it is also here possible to consider factorisations of
$C_{\C{Q}}$ in \refeq{eq:ev-cq-mult}, such as
\begin{equation}  \label{eq:fact-Q}
   C_{\C{Q}} = F F^* \approx F_a F_a^*, \quad \text{ with} \quad
    F, F_a :\C{E} \to \C{U}
\end{equation}
with some Hilbert space $\C{E}$, which lead again to factorisations of
\begin{equation}  \label{eq:fact-CQ}
   C = W^* W  \approx W_a^* W_a, \quad \text{ with } \quad
   W = F^*\Phi V^* \text{ and } 
   W_a = F_a^*\Phi V^* ,
\end{equation}
and representation on the space $\C{E}$; with the representing linear maps given by
$W^* = V \Phi^* F$ resp.\ $W_a^* = V \Phi^* F_a$.

Coming back to the situation where $C$ has a purely discrete spectrum
and a CONS of eigenvectors $\{v_m\}_m$ in $\C{U}$,
the map $B$ from the decomposition $C=B^* B$ can be used to
define a CONS $\{h_m\}_m$ in $\C{H}$:  $h_m := B C^{-1/2} v_m$, which is
 an eigenvector CONS of the operator $C_{\C{H}} := B B^*:\C{H}\to\C{H}$,
with $C_{\C{H}} h_m := \lambda_m h_m$, see   \citep{hgmRO-1-2018}.
From this follows a SVD of $B$ and $B^*$ in a manner analogous to \refeq{eq:XIV}.
The main result is  \citep{hgmRO-1-2018} that in the case of a nuclear
$C$  with necessarily purely discrete spectrum every factorisation leads
to a separated representation in terms of a series, and vice versa.  In
case $C$ is not nuclear, the representation of a ``parametric object''
via a linear map is actually more general \citep{hgmRO-1-2018, hgm-3-2018}
and allows to the rigorous and uniform treatment of also ``idealised''
objects, like for example Gaussian white noise on a Hilbert space.

In this instance of a discrete spectrum and a nuclear $C$ and hence nuclear $C_{\C{Q}}$,
the abstract equation $C_{\C{Q}} = \sum_k \lambda_k s_k \otimes s_k$ 
can be written in a more familiar form
in the case when the inner product on $\C{Q}$ is given
by a measure $\vpi$ on $\C{M}$.  It becomes for all $\vphi, \psi\in\C{Q}$:
\begin{multline*}  
   \ip{C_{\C{Q}}\vphi}{\psi}_{\C{Q}} = \sum_k \lambda_k \ip{\vphi}{s_k}_{\C{Q}}
   \ip{s_k}{\psi}_{\C{Q}} =    \ip{R^* \vphi}{R^* \psi}_{\C{U}} \\ 
   = \iint_{\C{M}\times\C{M}}  \vphi(\mu_1) \ip{r(\mu_1)}{r(\mu_2)}_{\C{U}}
    \psi(\mu_2)\; \vpi(\di \mu_1) \vpi(\di \mu_2) \\
   = \iint_{\C{M}\times\C{M}}  \vphi(\mu_1) \vkappa(\mu_1, \mu_2)
    \psi(\mu_2)\; \vpi(\di \mu_1) \vpi(\di \mu_2) \\
   = \iint_{\C{M}\times\C{M}}  \vphi(\mu_1) \left( \sum_k \lambda_k s_k(\mu_1)  s_k(\mu_2)
    \right) \psi(\mu_2)\; \vpi(\di \mu_1) \vpi(\di \mu_2).
\end{multline*}
This shows that $C_{\C{Q}}$ is really a Fredholm integral operator,
and its spectral decomposition is nothing but the familiar theorem 
of Mercer \citep{courant_hilbert} for the kernel
\begin{equation}  \label{eq:Mercer}
  \vkappa(\mu_1, \mu_2) = \sum_k \lambda_k s_k(\mu_1)  s_k(\mu_2) .
\end{equation}
Factorisations of $C_{\C{Q}}$ are then usually expressed as factorisations of the
kernel $\vkappa(\mu_1, \mu_2)$, which may involve integral transforms already envisioned
in \citep{Karhunen1947} --- see also the English translation \citep{Karhunen1947-e}:
\[
  \vkappa(\mu_1,\mu_2) = \int_{\C{Y}} \rho(\mu_1,y) \rho(\mu_2,y)\, \ops{n}(\di y), 
\]
where the ``factors'' $\rho(\mu, y)$ are measurable functions 
on the measure space $(\C{Y},\ops{n})$.  This is the classical analogue of
the general ``kernel theorem'' \citep{gelfand64-vol4}.

\paragraph{Connections to tensor products:}
Although not as obvious as for the case of a discrete spectrum in \refeq{eq:SVD-R-fdim},
\refeq{eq:spec-C-fdim}, and \refeq{eq:spec-CQ-fdim}; and \refeq{eq:XIII}, \refeq{eq:XIV},
such a connection is also possible in the general case of a non-discrete spectrum.
But as the spectral values in the continuous part have no corresponding eigenvectors,
one has to use the concept of generalised eigenvectors 
\citep{gelfand64-vol3, gelfand64-vol4, DautrayLions3, hgmRO-1-2018}.
Then it is possible to formulate the spectral theorem in the following way:
\begin{align}  \label{eq:spec-C-gen-ev}
   \ip{C u}{w}_{\C{U}} &= \int_{\D{R}^+} \lambda\, \ip{u}{v_\lambda} \ip{w}{v_\lambda}\, 
       \nu(\di \lambda) , \quad \text{ or in a weak sense}\\
       \label{eq:spec-C-gen-ev-wk}
    C &= \int_{\D{R}^+} \lambda\, (v_\lambda \otimes v_\lambda) \, \nu(\di \lambda),
\end{align}
with the spectral measure $\nu$ on $\D{R}^+$.
Observe the analogy, especially of \refeq{eq:spec-C-gen-ev-wk}, with \refeq{eq:XIII},
where the sum now has been generalised to an integral to account for the continuous spectrum.
\refeq{eq:spec-C-gen-ev} is for the case of a simple spectrum; in the more general
case of spectral multiplicity large than one,
the Hilbert space $\C{U}=\bigoplus_m \C{U}_m$ can be written as an orthogonal sum
\citep{gelfand64-vol3, gelfand64-vol4, DautrayLions3}
of Hilbert subspaces $\C{U}_m$, each invariant under the operator $C$, on which an
expression like  \refeq{eq:spec-C-gen-ev} holds, and on which the spectrum is simple.  
For the sake of brevity we shall
only consider the case of a simple spectrum now, and avoid writing the sums over $m$.
The difficulty in going from \refeq{eq:XIII} to \refeq{eq:spec-C-gen-ev-wk} is
that the values $\lambda$ in the truly continuous spectrum have no corresponding
eigenvector, i.e.\ $v_\lambda \notin \C{U}$, but it has to be found in a generally
larger space.
The possibility of writing an expression like \refeq{eq:spec-C-gen-ev} rests on the
concept of a ``rigged'' resp.\ ``equipped'' Hilbert space or Gel'fand triple.
This means that one can find \citep{gelfand64-vol4} a nuclear space 
$\C{K}\hookrightarrow \C{U}$ densely embedded in the Hilbert space $\C{U}$,
such that \refeq{eq:spec-C-gen-ev} holds for all $u, v \in \C{K}$.  This also
means that the generalised eigenvectors should be seen as linear
functionals on $\C{K}$.  As the subspace $\C{K}$ is densely embedded in $\C{U}$, 
it also holds that $\C{U}\hookrightarrow \C{K}^*$ is densely
embedded in the topological dual $\C{K}^*$ of $\C{K}$, i.e.\ one has the Gel'fand triple
\begin{equation}  \label{eq:Gelfand-triple}
   \C{K}\hookrightarrow \C{U} \hookrightarrow \C{K}^* .
\end{equation}
The generalised eigenvectors can now be seen as elements of the dual, 
$v_\lambda\in\C{K}^*$, where the generalised eigenvalue equation 
$C v_\lambda = \lambda v_\lambda$ holds after an appropriate extension of $C$.

If an expressions such as \refeq{eq:spec-C-gen-ev} or \refeq{eq:spec-C-gen-ev-wk}
have to be approximated numerically, it becomes necessary to evaluated the integral
in an approximate way.  The integral is really only over the spectrum of $\sigma(C)$
of $C$, as outside of $\sigma(C)$ the spectral measure $\nu$ vanishes.
Obviously, one would first split the spectrum $\sigma(C) = \sigma_d(C) \cup \sigma_c(C)$
into a discrete $\sigma_d(C)$ and a continuous part $\sigma_c(C)$.  On the discrete
part, the integral is just a sum as shown before.  On the continuous part,
the integral has to be evaluated by a quadrature formula.  Choosing quadrature
points $\lambda_z\in\sigma_c(C)$ and appropriate integration weights $w_z\in\D{R}$,
the integral can be approximated by
\[
  \int_{\sigma_c(C)} \lambda\, \ip{u}{v_\lambda} \ip{w}{v_\lambda}\, \nu(\di \lambda)
   \approx \sum_z w_z \lambda_z \ip{u}{v_{\lambda_z}} \ip{w}{v_{\lambda_z}},
\]
an expression very similar to the ones used in case of discrete spectra.

\paragraph{Further tensor products:}
Essentially, the constructions we have been investigating could be seen
as elements of the tensor product $\C{U}\otimes\C{Q}$, or extensions thereof
as in the preceding paragraph.  
Often one, or both of the spaces $\C{U}$ or $\C{Q}$, can be further so factored
in tensor products, say without loss of generality 
that  $\C{Q} = \C{Q}_I \otimes \C{Q}_{II}$.
This is for example the case for the white-noise modelling of random
fields \citep{kreeSoize86, Janson1997, boulder:2011, hgm-3-2018}, where
one has $\C{Q}= \bigotimes_{m=1}^\infty \C{S}_m$.  We just want to indicate
how this structure can be used for further approximation.

It essentially means that the whole foregoing is applied, instead on $\C{U}\otimes\C{Q}$,
on the tensor product $\C{Q}_I\otimes\C{Q}_{II}$.  Combined with the upper level
decomposition on $\C{U}\otimes\C{Q}$, one sees that this becomes one on
$\C{U}\otimes (\C{Q}_I\otimes\C{Q}_{II})$.
The bilinear forms \refeq{eq:IX-a} and \refeq{eq:IX} can thus be written
as tri-linear forms, making a direct connection to tensor products and multi-linear
forms \citep{Hackbusch_tensor}.  Like in the just cited example of random fields
\citep{boulder:2011, hgm-3-2018},
often this can be extended to higher order tensor products in a tree-like manner
--- by splitting $\C{U}$, or $\C{Q}_I$ resp. $\C{Q}_{II}$.  This leads to a
hierarchical structure encoded in this binary tree, with the top product $\C{U}\otimes\C{Q}$
the root of the tree, and the individual factors as leaves of the tree.
The higher the order of the tensor product, the better it is possible
to exploit dependencies in low-rank formats \citep{Hackbusch_tensor,
GrasedyckKressnerTobler2013}.  This has been recently also pointed out
in the tight connections between \emph{deep neural networks} 
\citep{CohenSha2016, KhrulkovEtal2018} and such
tensor decompositions, which come in different \emph{formats} 
or \emph{representations} \citep{Hackbusch_tensor}.
The indicated binary tree leads to what is known as a 
\emph{hierarchical Tucker}- or HT-format; but obviously the multi-factor
tensor product can be split also in a non-binary fashion, leading to more
general tree-based tensor formats \citep{falco-hackb-nouy-2018}.  A completely
flat tree structure with only root and leaves corresponds to the
well known \emph{canonic polyadic}- or CP-decomposition or format,
the original \emph{proper generalised decomposition} (PGD) falls into
this category \citep{chinestaPL2011, AmmarChinestaFalco2010, Falco:2012,
chinestaBook, chinestaWillcox2017}.

%
%
%
%
%
%


%

\section{Structure preservation}  \label{S:xmpls}
The foregoing development for a parametric model $r:\C{M}\to\C{U}$
did not assume anything more than that $\C{U}$ is a Hilbert space.
In \refS{parametric} it was already indicated on how to proceed if
$\C{U}$ is not a Hilbert space, but a more general topological vector
space.  The treatment so far preserves the linear structure of the
space $\C{U}$, and the approximations are using this linear structure as
well.  The tensor based representations using tensors of certain
rank already have a more difficult geometric structure \citep{falco-hackb-nouy-2015},
indeed a manifold structure \citep{billaud-falco-nouy-2017}.

But here the concern is about the structure of the image set of the parametric
object $r(\mu)$ and its preservation under the approximations or ROMs $r_a(\mu)$.
In case the image set--- here $\C{U}$---is not a vector space, but say a
differential manifold, things are bound to get more complicated; one
possible route of attack seems to use the previous linear methods like
the ones described here to map into the tangent spaces.  One instance
of this, which seems to be more accessible, is the case when the image set
is a Lie group $\C{G}$.  Then everything can be done in the tangent space
at the group identity, the Lie algebra $\F{g}$ of the Lie group $\C{G}$.  The
Lie algebra is a linear space, and one may take $\C{U}:=\F{g}$.  One
then has to map further from $\F{g}$ to $\C{G}$, but this can be achieved by
the canonical exponential map $\exp:\F{g}\to\C{G}$.  A representation or ROM then would
have the form 
\[
  \C{M}\xrightarrow{r, r_a} \C{U}=\F{g} \xrightarrow{\exp} \C{G}.
\]
This has the added advantage that interpolations along straight lines in
$\F{g}$, which like in any Euclidean or unitary space are also geodetics,
is mapped into interpolations along geodetics on the Riemannian manifold
structure on $\C{G}$.
We shall come back to a somewhat similar situation later in this section.

\paragraph{Vector fields:}
One of the probably simplest situations is when the image set has the structure
of $\C{V} = \C{U} \otimes \C{E}$, where
$\C{E}$ is a \emph{finite-dimensional} inner-product (Hilbert) space \citep{kreeSoize86}:
\begin{equation}   \label{eq:param-vect-r}
   \tns{r}:\C{M}\to\C{V} = \C{U} \otimes \C{E};\quad \tns{r}(\mu) = \sum_k r_k(\mu) \vek{r}_k,
\end{equation}
and the $r_k$ are maps $r_k(\mu)\C{M}\to\C{U}$ as before in \refS{parametric} and
\refS{correlat}, whereas the $\vek{r}_k$ are typically linearly independent vectors in $\C{E}$.
Often one wants to preserve the structure  $\C{V} = \C{U} \otimes \C{E}$;
one can think of this in the following way: $\C{U}$ is a space of scalar functions,
on some domain in Euclidean space, and $\C{E}$ are vectors from the associated vector
space. Hence one could call this a vector field in some sense.  The associated linear
map is then defined a bit differently, namely as
\[
   R_{\C{E}}: \C{U} \to \C{Q}\otimes\C{E}; \quad R_{\C{E}}: 
   u \mapsto \sum_k (R_k(\mu) u )\vek{r}_k,
\]
where the maps $R_k:\C{U}\to\C{Q}$ are defined as before in \refeq{eq:IX-a}.

The ``correlation'' can now be given by a bilinear form; namely
the densely defined map $C_{\C{E}}$ in $\C{V}=\C{U}\otimes\C{E}$ is defined
on elementary tensors $\tns{u} = u\otimes\vek{u}, \tns{v}=v\otimes\vek{v}
 \in \C{V}=\C{U}\otimes\C{E}$ as
\begin{equation}   \label{eq:vect-corr}
   \ip{C_{\C{E}}\tns{u}}{\tns{v}}_{\C{v}} :=   \sum_{k,j}  \ip{R_k(u)}{R_j(v)}_{\C{Q}}
    \,(\vek{u}^{\ops{T}}\vek{r}_k)\, (\vek{r}_j^{\ops{T}}\vek{v})
\end{equation}
and extended by linearity, where each $R_k:\C{U}\to\C{Q}$ is the map associated to $r_k(\mu)$
as before for just a single map $r(\mu)$.
It may be called the 
\emph{``vector correlation''}.  By construction it is self-adjoint and positive.
The corresponding kernel is not scalar, but has values in $\C{E}\otimes\C{E}$:
\begin{equation}   \label{eq:vect-kernel}
   \vek{\vkappa}_{\C{E}}(\mu_1,\mu_2) = \sum_{k,j} \ip{r_k(\mu_1)}{r_j(\mu_2)}_{\C{U}} 
   \; \vek{r}_k\otimes\vek{r}_j  .
\end{equation}
The eigenvalue problem on for an integral operator with such a kernel ---
representing the companion map --- is on $\C{W}=\C{Q}\otimes\C{E}$.

\paragraph{Coupled systems:}
A in some way similar situation is when the state space $\C{U}=\C{U}_1\times\C{U}_2$
comes from a combined or \emph{coupled} system \citep{hgmRO-2-2018}, and one wants
to conserve this information or structure.  The state is represented as
$\vek{u}=(u_1,u_2)$, and the natural inner product on such a
normal product space is 
\[ 
  \ip{\vek{u}}{\vek{v}}_{\C{U}} = \ip{(u_1,u_2)}{(v_1,v_2)}_{\C{U}} = \ip{u_1}{v_1}_{\C{U}_1}
           + \ip{u_2}{v_2}_{\C{U}_2}
\]
for $\vek{u}=(u_1,u_2), \vek{v}=(v_1,v_2) \in \C{U}$.  
This is for two coupled systems, labelled as
`1' and `2'.  The parametric map  is
\begin{equation}   \label{eq:param-coup-r}
   \vek{r}:\C{M}\to\C{U} = \C{U}_1 \times \C{U}_2;\quad 
     \vek{r}(\mu) = (r_1(\mu), r_2(\mu)).
\end{equation}
The associated linear map is
\begin{equation}   \label{eq:lin_map-coup}
   \vek{R}_c:\C{U}\to\C{Q}^2 = \C{Q} \times \C{Q};\quad 
     (\vek{R}_c(\vek{u}))(\mu) = (\ip{u_1}{r_1(\mu)}_{\C{U}_1}, \ip{u_2}{r_2(\mu)}_{\C{U}_2}).
\end{equation}
As before, these $\D{R}^2$ valued functions on $\C{M}$ are like two problem-adapted
co-ordinate systems on the joint parameter set, one for each sub-system.
From this one obtains the ``coupling correlation'' $\vek{C}_c$, 
again defined through a bilinear form
\begin{equation}   \label{eq:coup-corr}
   \ip{\vek{C}_{c} \vek{u}}{\vek{v}}_{\C{U}} :=   \sum_{j=1}^2  \ip{R_j(u_j)}{R_j(v_j)}_{\C{Q}} .
\end{equation}
The kernel is then a $2 \times 2$ matrix valued function in an integral operator
on $\C{W}=\C{Q}\times\C{Q}$:
\begin{equation}   \label{eq:coup-kernel}
   \vek{\vkappa}_{c}(\mu_1,\mu_2) =  
   \diag (\ip{r_k(\mu_1)}{r_k(\mu_2)}_{\C{U}_k}) .
\end{equation}

Other variations regarding  coupled systems are possible, see \citep{hgmRO-2-2018},
like when the parameter set $\C{M}=\C{M}_1\times\C{M}_2$ is a product.
The parametric map can then defined as
\begin{equation}   \label{eq:param-coup-12}
   \vek{r}:\C{M}=\C{M}_1\times\C{M}_2\to\C{U} = \C{U}_1 \times \C{U}_2;\quad 
     \vek{r}((\mu_1,\mu_2)) = (r_1(\mu_1), r_2(\mu_2)) ,
\end{equation}
with the associated linear map
\begin{equation}   \label{eq:lin_map-coup-12}
   \vek{R}:\C{U}\to\C{Q} = \C{Q}_1 \times \C{Q}_2;\quad 
     (\vek{R}(\vek{u}))(\mu) = (\ip{u_1}{r_1(\mu_1)}_{\C{U}_1}, \ip{u_2}{r_2(\mu_2)}_{\C{U}_2}).
\end{equation}
The correlation may be defined as before in \refeq{eq:coup-corr}, and also the
kernel on $\C{Q}=\C{Q}_1\times\C{Q}_2$ is as in \refeq{eq:coup-kernel},
but now the first diagonal entry is a function on $\C{M}_1 \times \C{M}_1$ only,
and analogous for the second diagonal entry.

\paragraph{Tensor fields:}
This is similar to the case of vector fields in that the state space is
$\C{W}=\C{U}\otimes \C{A}$, where $\C{U}$ is a space of scalar valued
functions on some set; and $\C{A} \subset \C{B}= \C{E}\otimes\C{E}$, where
$\C{E}$ is a finite-dimensional vector space \citep{hgmRO-1-2018}, and
$\C{A}$ is a manifold of tensors in the full tensor product $\C{B}$ of tensors
of even degree.  If $\C{A}$ were the full tensor product $\C{B}$, which is a linear
finite-dimensional space, there would be no difference to the case of vector fields.
But in case of tensors of even degree are often used in more special situations.
Obviously, such tensors may be identified with linear maps \citep{hgmRO-1-2018}
$\E{L}(\C{E}) \cong \C{B}$, which will be done here.  Therefore one may speak
of e.g.\ the manifold of \emph{special orthogonal} tensors, say $\ops{SO}(\C{E})$,
and of the manifold of \emph{symmetric positive definite} tensors $\ops{Sym}^+(\C{E})$.

We shall consider only these two mentioned examples.
The special orthogonal tensors are a Lie group $\C{A} := \ops{SO}(\C{E})$ with 
Lie algebra $\F{a} := \F{so}(\C{E})$, the \emph{skew-symmetric} tensors,
a free linear space.  For $\vek{S}\in\F{so}(\C{E})$, the exponential
map carries it onto $\exp(\vek{S})\in\ops{SO}(\C{E})$.
Therefore a parametric element in $\C{W}=\C{U}\otimes \C{A}$ can first
be represented as a parametric element in the linear space $\C{Z}=\C{U}\otimes \F{a}$,
where all the preceding statements on vector fields apply.  It is on this
intermediary representation that one can define ROMs.  Such a representation is then 
further mapped through exponentiation:
\[
   \exp_1 : \C{Z}=\C{U}\otimes \F{a} \ni u \otimes \vek{S} \mapsto
             u \otimes \exp(\vek{S}) \in \C{U}\otimes \C{A} = \C{W}.
\]
The associated linear map goes from $\C{Z}=\C{U}\otimes \F{a}$ to the linear space
$\C{Y}=\C{Q}\otimes \F{a}$; and again from here one would use an analogue of
the above exponential to map on $\C{Q}\otimes \C{A}$.

The positive definite tensors  $\C{A} := \ops{Sym}^+(\C{E})$ are not a classical
Lie algebra under multiplication, i.e.\ concatenation of linear maps, but rather
only a Riemannian manifold, geometrically a convex cone.  But there still is
an exponential map, carrying the free linear space of symmetric tensors 
$\F{a} := \F{sym}(\C{E})$ onto $\C{A} := \ops{Sym}^+(\C{E})$.  
In fact, for a $\vek{H}\in\F{sym}(\C{E})$, the exponential maps
it onto $\exp(\vek{H})\in\ops{Sym}^+(\C{E})$.  Thus we have recovered formally the same situation
as for orthogonal tensors just described, and the same procedures may be followed.

\section{Conclusion} \label{S:concl}
Parametric mappings $r:\C{M}\to\C{U}$ have been analysed with in a  variety of
settings via the associated linear map $R:\C{U}\to\C{Q}\subseteq\D{R}^{\C{M}}$.
It was shown that the associated linear map contains the full information
present in the parametric entity.  It is actually a mathematically more general
concept which allows one to define extreme or idealised such entities; this
is particularly relevant in the field of uncertainty quantification when one
has to deal with stochastic processes and random fields \citep{hgm-3-2018}.

So instead of analysing a parametric entity $r(\mu)$ and ist approximations
or ROMs $r_a(\mu)$ directly, one may take the cues on how to do this from
considering the associated linear maps $R$ and $R_a$.  One has to say
that in practical situations the associated linear maps are typically not simply
available explicitly, but they provide a conceptual framework on how to
deal with the situation.  And even though they are not directly available,
the desired quantities needed in such analyses are all in principle computable.

Very closely related to such an associated linear map $R$ is the so-called ``correlation
operator'' $C=R^*R$ and its companion $C_{\C{Q}} = R R^*$, both self-adjoint
and positive definite.  Their  spectral analysis turns out to be very helpful
in understanding the nature of such parametric entities, as well as possible
ROMs.  The very general nature and mathematical embedding of parametric entities, 
which also incorporates random fields, is shown in the different spectral
properties of the correlation operator.  Such generalised parametric entities
may yield correlation operators with continuous spectra --- as it typically
occurs for homogeneous random fields ---  and thus this needs
the full generality of the spectral analysis in rigged Hilbert spaces for
understanding the spectral analysis in terms of generalised tensor products.
Other factorisations of the correlation, such as $C = B^* B$, induce other
representations for the parametric entities, and any other representation or
re-parametrisation may be understood in these terms.

Preservation of certain structural properties is often very desirable.  Examples
are given to show how the general idea can be refined to reflect some linear
structures in the representation.  This even extends to non-linear manifolds
if they can be easily parametrised by linear spaces.  Lie groups with their
associated Lie algebras are one such example which is mentioned in a bit more detail.
This last point is especially relevant to the representation of spatially
varying or even random material properties, which are typically fields
of symmetric positive tensors.  A similar comment applies to ``orientation
fields'', which are spatially varying and possibly random fields of orthogonal
tensors.

Additionally it is explained how representations in tensor product spaces arise
naturally in such situations, and how this process can be cascaded to produce
a tree like structure for the analysis.  Low-rank tensor approximations can thus
be used as ROMs, and this certainly offers fresh new impulses.  The same applies
to machine learning and data-driven approaches, which obviously can also be
analysed with the proposed framework.  These deep learning methods have recently
been shown to be closely connected with low-rank tensor approximations, offering
some insights and avenues for their analysis.  With the proposed framework of
analysing such parametric entities via linear maps, we hope to introduce a fresh
point of view which may lead to new ideas on how to construct and analyse ROMs.

\ignore{           
\subsection{Test of fonts} \label{SS:testf}
\paragraph{Slanted:} --- serifs
\[ a, \alpha, A, \Phi, \phi, \vphi, \qquad 
     \vek{a}, \vek{\alpha}, \vek{A}, \vek{\Phi}, \vek{\phi}, \vek{\vphi} \]

\paragraph{Slanted:} --- sans serif
\[ \tns{a}, \tns{\alpha}, \tns{A}, \tns{\Phi}, \tns{\phi}, \tns{\vphi}, \qquad
  \tnb{a}, \tnb{\alpha}, \tnb{A}, \tnb{\Phi}, \tnb{\phi}, \tnb{\vphi} \]
   
\paragraph{Other:}  --- no small letters
\[ \C{E}, \C{Q}, \C{R}; \D{E}, \D{Q}, \D{R}; \E{E}, \E{Q}, \E{R}\]

\paragraph{Fraktur:}
\[ \F{e}, \F{q}, \F{r}; \F{E}, \F{Q}, \F{R} \]

\paragraph{Upright:} --- serifs, no small greek letters
\[ \mrm{a}, \mrm{A}, \mrm{\Phi},\qquad \mat{a}, \mat{A}, \mat{\Phi} \]

\paragraph{Upright:} --- sans serif, no small greek letters
\[ \ops{a}, \ops{A}, \ops{\Phi},\qquad \opb{a}, \opb{A}, \opb{\Phi} \]
}           

%
%
%
%
%
%





\providecommand{\bysame}{\leavevmode\hbox to3em{\hrulefill}\thinspace}
\providecommand{\MR}{\relax\ifhmode\unskip\space\fi MR }
\providecommand{\MRhref}[2]{%
  \href{http://www.ams.org/mathscinet-getitem?mr=#1}{#2}
}
\providecommand{\href}[2]{#2}



{ 
   \tiny
       \texttt{\RCSId} 
   }



\end{document}